\documentclass[a4paper,11pt]{article}

\usepackage{amsmath,amsfonts,amsthm,amssymb}
\usepackage[english]{babel}
\newcommand{\N}{\mathbb{N}}
\newcommand{\R}{\mathbb{R}}
\newcommand{\C}{\mathbb{C}}

\newcommand{\dimH}{\mathrm{dim}_\mathrm{H}\,}
\newtheorem{theorem}{Theorem}
\newtheorem{lemma}{Lemma}
\newtheorem{corollary}{Corollary}
\newtheorem{proposition}{Proposition}
\theoremstyle{definition}
\newtheorem{definition}{Definition}

\title{Large intersection classes on fractals}
\author{
David F\"arm\footnote{Centre for Mathematical Sciences, Lund University, Box \oldstylenums{118}, \textsc{se}-\oldstylenums{22100} Lund, Sweden, \texttt{david@maths.lth.se}}\ \ and
Tomas Persson\footnote{Institute of Mathematics, Polish Academy of Sciences, ulica \'Sniadeckich \oldstylenums{8}, Box \oldstylenums{21}, \oldstylenums{00}-\oldstylenums{956} Warszawa, Poland, \texttt{tomasp@impan.gov.pl}, Tomas Persson was supported by {\sc ec fp}\oldstylenums{6} Marie Curie ToK programmes {\sc spade}\oldstylenums{2}, and {\sc cody}.}
}

\begin{document}
\maketitle

\begin{abstract}
We consider limit sets of some conformal iterated function systems, and introduce classes of subsets of the limit set, with the property that the classes are closed under countable intersections and all sets in the classes have large Hausdorff dimension. We show some applications in ergodic theory, Diophantine approximation and complex dynamics.
\end{abstract}

\section{Introduction}

In \cite{falconer}, Falconer defined classes $\mathcal{G}^s$ of $G_\delta$-subsets of $\R^n$, for $0 < s \leq n$. Falconer proved that these classes are invariant under countable intersections, and that any set in $\mathcal{G}^s$ has at least Hausdorff dimension $s$. Moreover these classes are invariant under inverse images of bi-Lipschitz maps.

In this paper we will study classes of subsets of limit sets $\Lambda$ of conformal graph directed Markov systems. Simple examples of such sets $\Lambda$ are the interval, the middle third Cantor set, and Julia sets of polynomials of degree larger than 1. We are interested in such classes which have the intersection and dimension properties of $\mathcal{G}^s$.

As is easily seen, the complement of a set in $\mathcal{G}^s$ cannot contain an open non-empty set. Since the complement of $\Lambda$ may contain open non-empty sets, the classes of Falconer can not be directly used in our case. We will therefore introduce a modification of Falconer's classes, consisting of subsets of $\Lambda$ with the desired properties. This will be done in Section~\ref{sec:classdef}. Before this is done we need some notations and results, presented in Section~\ref{sec:gdms}.

In Section~\ref{sec:applications}, we consider some applications of these classes. We show that sets of points for which the ergodic averages of a H\"older continuous function $g$ does not converge, are in our classes of sets, and we will apply this to piecewise expanding interval maps and Julia sets of polynomials. We will also study subsets of the middle third Cantor set, with some Diophantine properties, and show that these sets are in our classes.

\section{Graph directed Markov systems} \label{sec:gdms}

Let $\Sigma_A$ be a transitive subshift of finite type over the alphabet $\{1, \ldots, q\}$ governed by a transition matrix $A$. We endow $\Sigma_A$ with the product topology.

Elements in $\Sigma_A$ will be denoted by $\boldsymbol{i} = i_1 i_2 \ldots$
The topology of $\Sigma_A$ is generated by the metric $d_\Sigma$ defined by $d_\Sigma (\boldsymbol{i}, \boldsymbol{j}) = 2^{-n}$, where $n$ is the smallest number such that $i_n \neq j_n$.

We define a cylinder to be a set of the form
\[
  C_{i_1 \ldots i_k} = \{\, \boldsymbol{j} \in \Sigma_A : j_1 \ldots j_k = i_1 \ldots i_k \,\}. 
\]
On $\Sigma_A$ we have the shift map $\sigma$ defined by $\sigma \colon i_1 i_2 \ldots \mapsto i_2 i_3 \ldots$

We will now define a {\em graph directed Markov system}. This can be done in a much more general way, see for instance \cite{Mauldin-Urbanski2}. However in this paper we will only consider the following simplification. Let $X$ be a compact and non-empty subset of $\R^n$, endowed with the usual metric. We assume that there are contractions on $X$ denoted by $(f_i)_{i=1}^q$ and $(f_{i,j})_{i,j}$, where $ij$ is a word in $\Sigma_A$. We will consider compositions of the contractions of the form $f_{i_1} \circ f_{i_1, i_2} \circ f_{i_{n-1}, i_n}$ where $i_1 i_2 \ldots$ is an element in $\Sigma_A$. We will call such a system of contractions together with $\Sigma_A$, a graph directed Markov system.

For any element $i_1 i_2 \ldots$ in $\Sigma_A$, there is a point $x \in X$ such that
\[
  \{ x \} = \bigcap_{n=1}^\infty f_{i_1} \circ f_{i_1, i_2} \circ \cdots \circ f_{i_{n-1}, i_n} (X).
\]
This defines a mapping $\pi \colon \Sigma_A \to X$. Since all $f_i$ and $f_{i,j}$ are contractions, $\pi$ is continuous. The limit set or attractor of the graph directed Markov system is the set $\pi (\Sigma_A)$ and will be denoted by $\Lambda_A$. It is compact since $\Sigma_A$ is compact and $\pi$ is continuous. We endow $\Lambda_A$ with the subset topology that comes from the topology of $X$.

We will now state some assumptions on the graph directed Markov system, that will be used in this paper.

First of all, we will only consider the case when all the contractions of the graph directed Markov system are conformal maps. Such a system is called a {\em conformal graph directed Markov system}.

The first two assumptions are
\begin{equation} \label{assumption:X}
\begin{minipage}{0.8 \textwidth}
The compact set $X$ is the closure of its interior, and the images of the interior of $X$ under the conformal contractions are disjoint.
\end{minipage}
\end{equation}
and
\begin{equation} \label{assumption:contraction}
\begin{minipage}{0.8 \textwidth}
We assume that there are numbers $0 < \lambda_1 < \lambda_2 < 1$ such that $\lambda_1 < \lVert \mathrm{d}_x f_{i,j} \rVert < \lambda_2$ and $\lambda_1 < \lVert \mathrm{d}_x f_i \rVert < \lambda_2$ for all $i,j$ and $x \in X$.
\end{minipage}
\end{equation}

We will need some bounded distortion property. Consider the following assumption. There is an $\alpha > 0$ such that
\begin{equation} \label{eq:distortion}
  \bigl| \lVert \mathrm{d}_x f_{i,j} \rVert - \lVert \mathrm{d}_y f_{i,j} \rVert \bigr| \leq \sup_{z \in X} \lVert (\mathrm{d}_z f_{i,j} )^{-1} \rVert^{-1} |x - y|^\alpha
\end{equation}
holds for all $i,j$ and $x,y$, and similarly for the maps $f_i$. In \cite{Mauldin-Urbanski2} it is shown that, under the other assumptions above, this assumption is satisfied if the dimension $n$ of the space is larger than one. We will therefore work with the following assumption.
\begin{equation} \label{assumption:distortion}
\begin{minipage}{0.8 \textwidth}
The dimension $n$ of the space is larger than one, or the dimension is one and there is an $\alpha > 0$ such that
\[
  \bigl| | f'_{i,j} (x) | - | f'_{i,j} (y) | \bigr| \leq \sup_{z \in X} | f'_{i,j} (z) | |x - y|^\alpha
\]
holds for all $i,j$ and $x,y$, and similarly for the maps $f_i$.
\end{minipage}
\end{equation}

For future use, we note that \eqref{eq:distortion} implies that the following bounded distortion property holds, see \cite{Mauldin-Urbanski2}.
There is a constant $\kappa$ such that for any sequence $\boldsymbol{i} \in \Sigma_A$, integer $n > 0$ and $x, y \in X$ it holds that
\begin{equation} \label{eq:distortion2}
  \kappa^{-1} \leq \frac{ \lVert \mathrm{d}_x (f_{i_1} \circ f_{i_1,i_2} \circ \cdots \circ f_{i_{n-1}, i_n}) \rVert }{ \lVert \mathrm{d}_y (f_{i_1} \circ f_{i_1,i_2} \circ \cdots \circ f_{i_{n-1}, i_n}) \rVert } \leq \kappa.
\end{equation}

\section{A class with large intersection} \label{sec:classdef}

We define a semi-metric $d$ on $\Sigma_A$ by $d (\boldsymbol{i}, \boldsymbol{j}) = | \pi (\boldsymbol{i}) - \pi (\boldsymbol{j}) |$, so that
\[
  d (C_{i_1 \ldots i_n}) = | \pi (C_{i_1 \ldots i_n}) |.
\]
Note however that this semi-metric does not necessarily generate the topology of $\Sigma_A$. In fact, this need not be a metric, since there can be different $\boldsymbol{i}$ and $\boldsymbol{j}$ such that $\pi (\boldsymbol{i}) = \pi (\boldsymbol{j})$.

We define outer measures $M_\infty^t$ on $\Sigma_A$ by
\[
  M_\infty^t (A) = \inf \biggl\{\, \sum_{i} d (C_i)^t : A \subset \bigcup_i C_i \, \biggr\},
\]
where each $C_i$ is a cylinder of some generation.

We will work mainly with $M_\infty^t$, but for some future use we also introduce the outer measures $M_\delta^t$ and $M^t$ on $\Sigma_A$, defined by
\[
  M_\delta^t (A) = \inf \biggl\{\, \sum_{i} d (C_i)^t : A \subset \bigcup_i C_i, \ d (C_i) \leq \delta \, \biggr\},
\]
where each $C_i$ is a cylinder of some generation, and $M^t (A) = \lim_{\delta \to 0} M_\delta^t (A)$. The $t$-dimensional Hausdorff measure is
\[
  H^t (A) = \lim_{\delta \to 0} \inf \biggl\{\, \sum_{i} d (U_i)^t : A \subset \bigcup_i U_i, \ d (U_i) \leq \delta \, \biggr\}.
\]
We note that $M^t (A) \geq M_\infty^t (A)$ for all sets $A$, and that $M^t$ is equivalent to the $t$-dimensional Hausdorff measure $H^t$.

The outer measures $M_\infty^t$, $M_\delta^t$, $M^t$ and $H^t$ are defined on subsets of $\Sigma_A$.
Using the projection $\pi \colon \Sigma_A \to I$, we project these outer measures to outer measures on $X$, and denote them  by $\mathfrak{M}_\infty^t$, $\mathfrak{M}_\delta^t$, $\mathfrak{M}^t$ and $\mathfrak{H}^t$.
Equivalently we may define
\[
  \mathfrak{M}_\infty^t (A) = \inf \Bigl\{\, \sum_i | \pi (C_i) |^t : A \subset \bigcup_i \pi (C_i) \, \Bigr\},
\]
where each $C_i$ is a cylinder, and similarly for the measures $\mathfrak{M}_\delta^t$, $\mathfrak{M}^t$ and $\mathfrak{H}^t$.

We will use the following classes of sets, which are modifications of Falconer's intersection classes from \cite{falconer}.

\begin{definition}
Let $\mathcal G^t (\Lambda_A)$, $0 < t \leq \dimH \Lambda_A$ be the class of $G_\delta$-sets $F \subset \Lambda_A$ such that
\[
  \mathfrak{M}^t_\infty ( F \cap \pi(C) ) = \mathfrak{M}_\infty^t ( \pi(C) )
\]
holds for all cylinders $C$.
\end{definition}

We are mainly interested in the classes $\mathcal G^t (\Lambda_A)$, but to develop the theory, we prefer to work with classes of subsets of $\Sigma_A$.
Therefore we also introduce the following classes.

\begin{definition}
Let $\mathcal G^t (\Sigma_A)$, $0 < t \leq \dimH \Lambda_A$ be the class of $G_\delta$-sets $F \subset \Sigma_A$ such that
\[
  M^t_\infty (F \cap C) = M_\infty^t (C)
\]
holds for all cylinders $C$.
\end{definition}

Our main result on these classes is the following theorem.

\begin{theorem} \label{the:classes}
  The classes $\mathcal{G}^t (\Lambda_A)$ and $\mathcal{G}^t (\Sigma_A)$ are closed under countable intersections and the Hausdorff dimension of any set in one of these classes is at least $t$.
\end{theorem}

The proof of Theorem~\ref{the:classes} is in Section~\ref{sec:classproof}.

\section{Applications} \label{sec:applications}

\subsection{Non-typical points in ergodic theory}

In this section we will state some results that are generalisations of the results in \cite{farm} and \cite{schmeling2}.


For functions $g \colon \Sigma_A \to \R$, consider the ergodic averages
\begin{align}\label{average}
  \frac{1}{n} \sum_{k=0}^{n-1} g (\sigma^k \boldsymbol i),
\end{align}
as $n\to \infty$. For $x\in \Lambda_A$ such that $\pi(\boldsymbol i )=x$, let $A_g(x)$ be the set of accumulation points for the expression (\ref{average}).

Let $\Sigma$ be a shift space. We will say that a function $g \colon \Sigma \to \R$ is H\"older continuous if there are constants $K$ and $\alpha > 0$ such that $|g (\boldsymbol{i}) - g (\boldsymbol{j})| < K d_\Sigma (\boldsymbol{i}, \boldsymbol{j})^\alpha$ holds for all $\boldsymbol{i}, \boldsymbol{j} \in \Sigma$. Denote by $C_\alpha^K (\Sigma)$ the family
\[
  C_\alpha^K (\Sigma) = \{\, g \colon \Sigma \to \R : |g (\boldsymbol{i}) - g (\boldsymbol{j})| < K d_\Sigma (\boldsymbol{i}, \boldsymbol{j})^\alpha,\ \forall \boldsymbol{i}, \boldsymbol{j} \in \Sigma \,\},
\]
and by $C_\alpha^K (X)$ the family
\[
  C_\alpha^K (X) = \{\, \phi \colon X \to \R : |\phi (x) - g (y)| < K |x - y|^\alpha,\ \forall x, y \in X \,\}.
\]
We note that if $\phi \colon X \to \R$ is in $C_\alpha^K (X)$, then $\phi \circ \pi^{-1}$ is in $C_{\alpha_0}^{K_0} (\Sigma)$, for some constant $\alpha_0$ and $K_0$ that do not depend on $\phi$, \mbox{i.e.} $\pi^{-1}$ transfers functions of $C_\alpha^K (X)$ to functions of $C_{\alpha_0}^{K_0} (\Sigma)$.

We will prove the following.

\begin{theorem}\label{the:maintheorem}
  Consider a conformal graph directed Markov system, satisfying the assumptions \eqref{assumption:X}, \eqref{assumption:contraction} and \eqref{assumption:distortion}. 
  For any sequence $(g_i)_{i=1}^\infty$ of real valued continuous functions on $\Sigma_A$, and each sequence $(x_i)_{i=1}^\infty$ of points in $\R$, it holds that
  \[
    \dimH ( \cap_{i=1}^\infty \{ \, y \in \Lambda_A : x_i \in A_{g_i}(y) \, \} )
    = \inf_i \dimH ( \{ \, y \in \Lambda_A : x_i \in A_{g_i}(y) \, \} ).
  \]
\end{theorem}


We will combine Theorem \ref{the:maintheorem} with Theorem~9 of Barreira's and Saussol's article \cite{barreirasaussol}, to get the following corollary.

\begin{corollary} \label{cor:fulldim}
  Consider a conformal graph directed Markov system, satisfying the assumptions \eqref{assumption:X}, \eqref{assumption:contraction} and \eqref{assumption:distortion}.
  Let $(g_i)_{i=1}^\infty$ be a sequence of H\"older continuous maps on $\Sigma_A$ such that no $g_i$ is cohomologous to a constant. Then the set of points for which the ergodic averages do not converge for any of the maps $g_i$ has the same Hausdorff dimension as $\Lambda_A$.
\end{corollary}

Theorem \ref{the:maintheorem} and Corollary \ref{cor:fulldim} are generalisations of results in \cite{farm}, where the case $\Lambda_A=[0,1]$ is considered, \mbox{i.e.} where there are no holes. Moreover, in \cite{farm}, each function $g$ is the characteristic functions of a cylinder $C_{w_1\dots w_k}$, so that the ergodic average (\ref{average}) is the frequency of occurrence of the word $w_1\dots w_k$ in a sequence $\boldsymbol{i} \in \Sigma_A$. 

To prove Theorem~\ref{the:maintheorem} it is sufficient to prove that the sets $\{ \, y \in \Lambda_A : x_i \in A_{g_i}(y) \, \}$ are in $\mathcal{G}^t (\Lambda_A)$ for any $t$ smaller than the dimension of the set.
The proofs of Theorem~\ref{the:maintheorem} and Corollary~\ref{cor:fulldim} are in Section~\ref{sec:proof}.

\begin{corollary} \label{cor:piecewisemonotone}
  Let $0 = a_0 < \cdots < a_n = 1$, and let $f$ be monotone and $C^2$ on each of the intervals $(a_k, a_{k+1})$, with $\inf |f'| > 1$ and $\sup |f'| < \infty$. Let $(g_i)_{i=1}^\infty$ be H\"older continuous functions on $[0, 1]$, that are not cohomologous to a constant. Assume that there are constants $K$, $V > 0$ and $\alpha > 0$ such that all $g_i$ are in $C_\alpha^K ([0,1])$ and $\sup g_i -  \inf g_i > V$.
Then the set of points in $[0, 1]$, such that no ergodic average of any $g_i$ converges, has Hausdorff dimension $1$.
\end{corollary}

In \cite{schmeling2} another method is used to prove theorems similar to Corollary~\ref{cor:piecewisemonotone}, but only for finitely many maps $g_i$ at a time. Note also that Corollary~\ref{cor:juliaset} below, is in the same spirit.

Corollary~\ref{cor:piecewisemonotone} also generalises Corollary~3 in \cite{FarmPerssonSchmeling}, by taking much more general maps into consideration. However, the result in \cite{FarmPerssonSchmeling} is stronger in that it actually proves that the countable intersection is in $\mathcal{G}^1([0,1])$, and so, the map $f$ does not need to be fixed.

Let us outline the proof of Corollary~\ref{cor:piecewisemonotone}.

\begin{proof}[Proof of Corollary~\ref{cor:piecewisemonotone}]
  The system $(f, I)$ generates a subshift $\Sigma_f$ on an alphabet of $n$ symbols. Given $\varepsilon > 0$, by Proposition~1 in \cite{HofbauerRaithSimon}, we can approximate $\Sigma_f$ by a subshift of finite type $\Sigma_A$, such that the set $K \subset [0, 1]$, corresponding to $\Sigma_A$, have Hausdorff dimension at least $1 - \varepsilon$. The construction of $\Sigma_A$ is by forbidding some cylinders of generation $n$.

The functions $g_i \circ \pi^{-1}$ are all in $C_{\alpha_0}^{K_0} (\Sigma_f)$, for some constant $K_0$ and $\alpha_0 > 0$.
We need to show that $\Sigma_A$ can be choosen such that no $g_i$ restricted to $\Sigma_A$ is cohomologous to a constant. To conclude this, it is sufficient to show that $\Sigma_A$ can be choosen such that no $g_i$ is constant on $\Sigma_A$. To get a contradiction, assume that for each $\Sigma_A$ there is a $g$ satisfying the assumptions of the theorem, and such that $g$ is constant on $\Sigma_A$. Then there is a point $\boldsymbol{i} \in \Sigma_f \setminus \Sigma_A$ such that $|g \circ \pi^{-1} (\boldsymbol{i}) - g \circ \pi^{-1} (\boldsymbol{j})| \geq V/2$ for all $\boldsymbol{j} \in \Sigma_A$. Since $\Sigma_A$ is constructed by forbidding cylinders of length $n$, we can choose $\boldsymbol{j}$ so that $d_\Sigma (\boldsymbol{i}, \boldsymbol{j}) \leq 2^{-n+1}$. Let $n \to \infty$ to get a contradiction to $g \circ \pi^{-1} \in C_{\alpha_0}^{K_0} (\Sigma_f)$.

Hence, we can choose $\Sigma_A$ such that no $g_i$, restricted to $\Sigma_A$, is cohomologous to a constant.
Now Corollary~\ref{cor:piecewisemonotone} follows by Corollary~\ref{cor:fulldim}, and letting $\varepsilon \to 0$.
\end{proof}

\subsection{Diophantine approximation and Cantor sets}

Let $K$ be the middle third Cantor set.
With $n= 3$, $\Sigma_A = \{0, 2\}^\mathbb{N}$, and an appropriate definition of the iterated function system, we have that $K = \Lambda_A$.
Let $\alpha > 1$ and consider the sets
\[
  W (\alpha) = \Bigl\{\, x \in [0, 1] : \Bigl| x - \frac{p}{q} \Bigr| < q^{-\alpha} \text{ for infinitely many } p \in \mathbb{Z}, q = 3^k \,\Bigr\}.
\]
These sets were studied in the paper \cite{Levesley-Salp-Velani} by Levesley, Salp and Velani.
They proved that the Hausdorff dimension of $W (\alpha) \cap K$ is $\frac{1}{\alpha} \frac{\log 2}{\log 3}$.
We will strengthen this result and prove the following theorem.

\begin{theorem} \label{the:diophantine}
  $W (\alpha) \cap K$ is in $\mathcal{G}^t (K)$ for all $t < \frac{1}{\alpha} \frac{\log 2}{\log 3}$.
\end{theorem}

The proof of this theorem is in Section~\ref{sec:diophantineproof}.

As a corollary of Theorem~\ref{the:maintheorem} and \ref{the:diophantine} we get that the set of points in $K \cap W (\alpha)$ such that the frequency of any word is undefined, has Hausdorff dimension $\frac{1}{\alpha} \frac{\log 2}{\log 3}$.

\subsection{Julia sets}

In this section we will apply Corollary~\ref{cor:fulldim} to non-typical points of the Julia set $J (f)$ of a polynomial $f (z)$.

\begin{corollary} \label{cor:juliaset}
  Let $f (z)$ be a polynomial of degree larger than $1$, such that $|f'| > 1$ on the Julia set $J (f)$. Let $(g_k)_{k=1}^\infty$ be a sequence of H\"older continuous functions $g_k \colon J(f) \to \C$, such that no $g_k$ is cohomologous to a constant with respect to $(f, J(f))$. Then the set of points in $J (f)$ such that the ergodic averages of $g_k$ does not converge for any $g_k$, has the same Hausdorff dimension as $J (f)$.
\end{corollary}

\begin{proof}
  Since $f$ is of degree $d > 1$, the Julia set is compact and non-empty. There are only finitely many critical points since $f$ is a polynomial. Since $J (f)$ is compact and invariant, and $J (f)$ contains no critical point since $f' \ne 0$ on $J (f)$, there exists an $R > 0$ such that the distance from $J (f)$ to any image of a critical point is a least $R$, and $| f' | > 1$ on any ball $B_R (x)$, $x \in J(f)$.
  Moreover, the compactness of $J (f)$ implies that there exists a constant $c$ such that $|f' (z)| > c > 1$ for all $z \in J(f)$. We can even choose $c > 1$ and $0 < r < R$ such that if $x \in J(f)$, then the pre-image of the ball $B_r (x)$ of radius $r$ around $x$ satisfies
  \[
    f^{-1} (B_r (x)) \subset \bigcup_{y \in f^{-1} (x)} B_{r/c} (y),
  \]
  and $| f' | > c$ on $B_r (x)$.

  We let $U$ be an open neighborhood around $J (f)$ defined by
  \[
    U = \bigcup_{x \in J (f)} B_r (x),
  \]
  and put $X = \overline{U}$. 

  Now $X$ contains no image of a critical point, so we can define inverse branches $(f_j)_{j=1}^d$ of $f$ on $X$, and put $f_{i,j} = f_j$. It is clear that $f_j$ are conformal contractions on $X$. Moreover $f_j \colon X \to X$. Indeed we have
  \[
    f^{-1} (X) = \overline{f^{-1} (U)} \subset \overline{ \bigcup_{x \in J (f)} B_{r/c} (x) } \subset \overline{U} = X,
  \]
  so, $f_j (X) \subset X$.

  Note that some $f_j$ may be discontinuous at the curves $f ( \partial f_j (X))$. However $f ( \partial f_j (X))$ are piecewise smooth curves and if necessary, we can cut up $X$ along these curves, so that $f_j$ satisfies the assumptions \eqref{assumption:X}, \eqref{assumption:contraction} and \eqref{assumption:distortion}.
  Now, Corollary~\ref{cor:fulldim} finishes the proof.
\end{proof}

\section{Proof of Theorem~\ref{the:classes}} \label{sec:classproof}

Let us first comment on the relation between the classes $\mathcal G^t (\Lambda_A)$ and $\mathcal G^t (\Sigma_A)$.
If $U \subset X$ is an open set, then $\pi^{-1} (U) \subset \Sigma_A$ is an open set.
Hence, if $F \subset X$ is a $G_\delta$-set, then so is $\pi^{-1} (F) \subset \Sigma_A$, since $\pi$ is continuous.
We are going to show that the class $\mathcal G^t (\Sigma_A)$ is closed under countable intersections.
Since countable intersections of $G_\delta$-sets are $G_\delta$-sets, this implies that the class $\mathcal G^t (\Lambda_A)$ is also closed under countable intersections.
Hence, to prove Theorem~\ref{the:classes}, we only need to prove the intersection statement for the class $\mathcal{G}^t (\Sigma_A)$. This is also true for the statement on the Hausdorff dimension.

Define the potential $\phi \colon \Sigma_A \to \R$ by $\phi \colon \boldsymbol{i} = (i_k)_{k=1}^\infty \mapsto \log \lVert \mathrm{d}_{\pi (\boldsymbol{i})} f_{i_1, i_2} \rVert$.
We also let
\[
  S_n \phi ( \boldsymbol{i} ) = \log \lVert \mathrm{d}_{\pi (\boldsymbol{i})} f_{i_1} \rVert + \sum_{k=1}^{n-1} \phi (\sigma^k (\boldsymbol{i})).
\]
Then
\[
  e^{\inf_{\boldsymbol{i} \in C_{j_1 \ldots j_n}} S_n \phi (\boldsymbol{i})} \leq d (C_{j_1 \ldots j_n}) \leq e^{\sup_{\boldsymbol{i} \in C_{j_1 \ldots j_n}} S_n \phi (\boldsymbol{i})}.
\]

Since $\phi$ is a H\"older continuous potential, we can write the pressure of $s \phi$ as
\begin{align*}
  P (s \phi)
  &= \lim_{n \to \infty} \frac{1}{n} \log \biggl( \sum_{C_{j_1 \ldots j_n} } e^{\inf_{\boldsymbol{i} \in C_{j_1 \ldots j_n} } s S_n \phi (\boldsymbol{i})} \biggr) \\
  &= \lim_{n \to \infty} \frac{1}{n} \log \biggl( \sum_{C_{j_1 \ldots j_n} } e^{\sup_{\boldsymbol{i} \in C_{j_1 \ldots j_n} } s S_n \phi (\boldsymbol{i})} \biggr),
\end{align*}
where the sums run over all cylinders $C_{j_1 \ldots j_n}$ of generation $n$, see \cite{Mauldin-Urbanski} and \cite{Mauldin-Urbanski2}.
It is known, see for instance \cite{Mauldin-Urbanski2}, that $P (s \phi) = 0$ is equivalent to $\dimH (\Lambda_A) = s$.

We note that for each non-empty cylinder $C$ we have
\begin{align*}
  P (s \phi)
  & = \lim_{n \to \infty} \frac{1}{n} \log \biggl( \sum_{C_{j_1 \ldots j_n} \subset C} e^{\inf_{\boldsymbol{i} \in C_{j_1 \ldots j_n} } s S_n \phi (\boldsymbol{i})} \biggr) \\
  & = \lim_{n \to \infty} \frac{1}{n} \log \biggl( \sum_{C_{j_1 \ldots j_n} \subset C} e^{\sup_{\boldsymbol{i} \in C_{j_1 \ldots j_n} } s S_n \phi (\boldsymbol{i})} \biggr),
\end{align*}
so given $s < \dimH (\Lambda_A)$, we have $P (s \phi) > 0$ and there are constants $m = m(s)$ and $c > 0$ such that
\begin{equation} \label{eq:positive}
  \biggl( \sum_{C_{j_1 \ldots j_n} \subset C_a} e^{\inf_{\boldsymbol{i} \in C_{j_1 \ldots j_n} }
  s S_n \phi (\boldsymbol{i})} \biggr) > e^{c n}> 1,
\end{equation}
for all cylinders $C_a$ and $n \geq m$.

We observe that if $C$ is a cylinder, then it need not be true that $M_\infty^t (C) = d (C)^t$, since $C$ might not be the optimal cover of $C$. However, because of \eqref{eq:positive}, we only need to consider covers of $C$ with cylinders that are at most of $m = m(t)$ generations higher than $C$. Since there are only finitely many such covers, we conclude that there is a largest constant $0 < c_t \leq 1$ such that
\begin{equation} \label{eq:measureofcylinder}
  M_\infty^t (C) \geq c_t d (C)^t
\end{equation}
holds for all cylinders $C$. A similar statement is of course true for the measure $\mathfrak{M}_\infty^t$. Note also that we always have $M_\infty^t (C) \leq d (C)^t$.

The proof of Theorem~\ref{the:classes}, will be divided into a sequence of Lemmata.
The proofs of these Lemmata are in most cases quite similar to the corresponding ones in Falconer's paper \cite{falconer}.
The following lemma lets us extend the property that defines $\mathcal G^t (\Lambda_A)$ to any open set.

\begin{lemma} \label{lem:openset}
  If\/ $c > 0$ and $F$ is a set such that
  \[
    M_\infty^t (F \cap C) \geq c M_\infty^t (C)
  \]
  holds for all cylinders $C$, then
  \[
    M_\infty^t (F \cap U) \geq c M_\infty^t (U)
  \]
  holds for all open sets $U$.
\end{lemma}

\begin{proof}
  Let $U$ be open. Then we can write $U$ as a disjoint union $U = \bigcup_{i}^\infty C_i$ where each $C_i$ is a cylinder.
  Let $(D_j)$ be cylinders covering $F \cap U$. We can assume that this cover is disjoint.

  By the net property of cylinders, for any $i$, either there are no $D_j \subset C_i$ and instead $C_i \subset D_j$ for some $j$,
  or there are $D_j \subset C_i$ and
  \[
    \bigcup_{D_j \cap C_i \ne \emptyset} D_j \subset C_i.
  \]
  In case there are $D_j \subset C_i$ we have that
  \begin{equation} \label{eq:subsetsofC}
    \sum_{D_j \subset C_i} d (D_j)^t
    \geq M_\infty^t (F \cap C_i)
    \geq c M_\infty^t (C_i).
  \end{equation}

  If we take $(D_j)$ and for each $i$ such that there are $D_j\subset C_i$, we replace all these $D_j$ by one copy of $C_i$, then we get a disjoint cover $(E_r)$ of $U$, where all $E_r$ are cylinders.
  Indeed, any $C_i$ from $U = \bigcup_{i}^\infty C_i$ is either a subset of some element $E_r = D_j$, or it is itself an element $E_r$. 

  Using \eqref{eq:subsetsofC} and the fact that $c \leq 1$, we get
  \[
    \sum d (D_j)^T \geq c \sum M_\infty^t (E_r) \geq c M_\infty^t (U).
  \]
  Since $(D_j)$ is arbitrary, this proves the lemma.
\end{proof}

We will also need the following lemma.

\begin{lemma} \label{lem:removeconstant}
  Let $F$ be a set such that there is a constant $c > 0$ such that
  \[
    M_\infty^{t_0} (F \cap C) \geq c M_\infty^{t_0} (C)
  \]
  holds for all cylinders $C$. Then
  \[
    M_\infty^t (F \cap C) \geq M_\infty^t (C) 
  \]
  holds for all cylinders $C$, and $0 < t \leq t_0$. 
\end{lemma}

\begin{proof}
  We start by proving that $M_\infty^t (F \cap C) \geq M_\infty^t (C)$ holds for $t = t_0$. 
  We then use this result to obtain $M_\infty^t (F \cap C) \geq M_\infty^t (C)$ when $t < t_0$.

  Let $(C_i)_i$ be a collection of cylinders covering $F \cap C$.
  We may assume that the sets $C_i$ are pairwise disjoint.
  Since $M_\infty^{t_0} (F \cap C)$ is finite,
  we may assume that $\sum d (C_i)^{t_0}$ is finite.
  We let $I(m) = \{\, i : C_i \text{ is of generation } m \,\}$.
  Then, for any $\varepsilon > 0$, there is an $m_0$ such that
  \begin{equation} \label{eq:m0partition}
    \sum_{m \geq m_0} \sum_{i \in I (m)} d (C_i)^{t_0}
    = \sum_i d (C_i)^{t_0} - \sum_{m < m_0} \sum_{i \in I (m)} d (C_i)^{t_0}
    < \varepsilon.
  \end{equation}

  We let $(D_j)$ denote a finite cover of the cylinder $C$ by cylinders,
  such that for any $j$ holds either
  \begin{itemize}
    \item[i)] $D_j = C_i$ for some $i$ and $D_j$ is of generation less
    than $m_0$.
  \end{itemize}
  or
  \begin{itemize}
    \item[ii)] $D_j$ is of generation $m_0$, and those $C_i$ that intersect
    $D_j \cap F$ are contained in $D_j$.
  \end{itemize}

  Let $Q(j) = \{\, i : C_i \subset D_j \,\}$. If $j$ satisfies i) then
  \begin{equation} \label{eq:ifi}
    \sum_{i \in Q(j)} d (C_i)^{t_0} = d (D_j)^{t_0}.
  \end{equation}

  For those $j$ that satisfies ii), we have for $i \in Q(j)$ that
  \begin{equation} \label{eq:ifii}
    \sum_{i \in Q(j)} d (C_i)^{t_0}
    \geq M_\infty^{t_0} (F \cap D_j)
    \geq c M_\infty^{t_0} (D_j)
    \geq c c_{t_0} d (D_j)^{t_0}.
  \end{equation}

  Using \eqref{eq:m0partition}, \eqref{eq:ifi} and \eqref{eq:ifii}, we conclude that
  \begin{align*}
    \sum_i d (C_i)^{t_0}
    &= \sum_{m < m_0} \sum_{i \in I (m)} d (C_i)^{t_0} + \sum_{m \geq m_0} \sum_{i \in I (m)} d (C_i)^{t_0} \\
    &= \sum_{m < m_0} \sum_{i \in I (m)} d (C_i)^{t_0} + c^{-1} c_{t_0}^{-1} \sum_{m \geq m_0} \sum_{i \in I (m)} d (C_i)^{t_0} \\
    & \hspace{3cm} + (1 - c^{-1} c_{t_0}^{-1}) \sum_{m \geq m_0} \sum_{i \in I (m)} d (C_i)^{t_0} \\
    &\geq \sum_j d (D_j)^{t_0} + (1 - c^{-1} c_{t_0}^{-1}) \varepsilon
    \geq M_\infty^{t_0} (C) +  (1 - c^{-1} c_{t_0}^{-1}) \varepsilon.
  \end{align*}
  Hence, $M_\infty^{t_0} (C \cap F) \geq M_\infty^{t_0} (C) +  (1 - c^{-1} c_{t_0}^{-1}) \varepsilon$.
  Since $\varepsilon$ is arbitrary, we conclude that $M_\infty^{t_0} (C \cap F) \geq M_\infty^{t_0} (C)$.

  We now turn to the case $0 < t < t_0$.
  Using what we proved above, we note that it is sufficient to prove that there is a constant $c'$ such that $M_\infty^t (C \cap F) \geq c' M_\infty^t (C)$ holds for all cylinders $C$.

  Let $C$ be an $n$-cylinder. Take $n_0 \geq n$ such that $\lambda_2^{n_0 + 1} \leq \lambda_1^{n + 1}$, or equivalently
  \begin{equation} \label{eq:n0}
    \lambda_2^{(n_0 + 1)(t - t_0)} \geq \lambda_1^{(n + 1)(t - t_0)}.
  \end{equation}
  Let $(C_i)_i$ be a collection of cylinders covering $F \cap C$.
  We may assume that the sets $C_i$ are pairwise disjoint.
  We let $(D_j)$ denote a finite cover of the cylinder $C$ by cylinders, such that for any $j$ holds either
  \begin{itemize}
    \item[i)] $D_j = C_i$ for some $i$ and $D_j$ is of generation less
    than $n_0$.
  \end{itemize}
  or
  \begin{itemize}
    \item[ii)] $D_j$ is of generation $n_0$, and those $C_i$ that intersect
    $D_j \cap F$ are contained in $D_j$.
  \end{itemize}

  Let $Q(j) = \{\, i : C_i \subset D_j \,\}$. If $j$ satisfies i) then
  \[
    \sum_{i \in Q(j)} d (C_i)^t = d (D_j)^t = d (D_j)^{t - t_0} d (D_j)^{t_0} \geq d (C)^{t - t_0} d (D_j)^{t_0}.
  \]

  If $j$ satisfies ii) then, for $i \in Q(j)$, holds
  \begin{align*}
    d (C_i)^t &= d (C_i)^{t - t_0} d (C_i)^{t_0}
    \geq (\lambda_2^{n_0 + 1})^{t - t_0} d (C_i)^{t_0} \\
    &\geq (\lambda_1^{n + 1})^{t - t_0} d (C_i)^{t_0}
    \geq d (C)^{t - t_0} d (C_i)^{t_0}, 
  \end{align*}
  by \eqref{eq:n0}. Hence, if $j$ satisfies ii), then
  \begin{align*}
    \sum_{i \in Q(j)} d (C_i)^t
    &\geq d (C)^{t - t_0} \sum_{i \in Q(j)} ^{t_0} (C_i)^{t_0}
    \geq d (C)^{t - t_0} M_\infty^{t_0} (F \cap D_j) \\
    &\geq d (C)^{t - t_0} M_\infty^{t_0} (D_j)
    \geq c_{t_0} d (C)^{t - t_0} d (D_j)^{t_0}.
  \end{align*}
  It follows that
  \begin{align*}
    \sum_i d(C_i)^t
    &\geq c_{t_0} d (C)^{t - t_0} \sum_j d (D_j)^{t_0}
    \geq c_{t_0} d (C)^{t - t_0} M_\infty^{t_0} (C) \\
    &\geq c_{t_0}^2 d (C)^t
    \geq c_{t_0}^2 M_\infty^{t_0} (C)
  \end{align*}
  This proves that
  \[
    M_\infty^t (F \cap C) \geq c_{t_0}^2 M_\infty^t (C) 
  \]
  holds for all cylinders $C$, and $0 < t < t_0$.
\end{proof}

The proof of Theorem~\ref{the:classes} is finished, if we prove the following two propositions.

\begin{proposition} \label{pro:dimension}
  If $F \in \mathcal G^t (\Sigma_A)$, then $\dimH (F) \geq t$ with respect to the metric $d$ on $\Sigma_A$ and thereby $\dimH (\pi(F)) \geq t$ with respect to Euclidean metric on $[0,1)$.
\end{proposition}

\begin{proposition} \label{pro:intersectionproperty}
  If $F_i \in \mathcal G^t (\Sigma_A)$ for all $i\in \N$, then
  \[
    M_\infty^t ( \bigcap_{i=1}^\infty F_i \cap U ) = M_\infty^t (U),
  \]
  for all open $U$, and thereby $\bigcap_{i=1}^\infty F_i  \in \mathcal G^t (\Sigma_A)$.
\end{proposition}

\begin{proof}[Proof of Proposition~\ref{pro:dimension}]
  Clearly $M_\infty^t (F) > 0$. Hence $M^t (F) > 0$, and since $M^t$ is equivalent to the $t$-dimensional Hausdorff-measure, we conclude that $H^t (F) > 0$ and $\dimH (F) > t$.
\end{proof}

\begin{proof}[Proof of Proposition~\ref{pro:intersectionproperty}]
  The proof of this proposition is very similar to the corresponding proof in
  \cite{falconer}. It is based on an increasing set lemma from Roger's book \cite{rogers}.

  Let $E\subset \Sigma_A$ and $\delta > 0$. We will denote by $E_{(-\delta)}$ the set
  \[
    E_{(-\delta)} = \{\, x \in E : \inf_{y \not \in E} d (x, y) > \delta \,\},
  \]
  and note that $E_{(-\delta)}$ is an open set. 

  Take $\varepsilon > 0$.
  Let first $\{F_i\}_{i=1}^\infty$ be a decreasing sequence of open sets,
  with the property that
  \[
    M_\infty^t (F_i \cap U) \geq M_\infty^t (U)
  \]
  holds for any open set $U$.
  We will choose a sequence of numbers $(\delta_k)$ and open sets $(U_k)$, such that
  \begin{align*}
    U_0 &= U, \\
    U_k &= (F_k \cap U_{k-1})_{(-\delta_k)},
  \end{align*}
  and $M_\infty^t (U_k) > M_\infty^t (U) - \varepsilon$.

  We will choose $\delta_k$ and $U_k$ inductively. Assume that they have been chosen
  for $k = 1, \ldots, n$. Since $F_n\in \mathcal G^t (\Lambda_A)$, Lemma \ref{lem:openset} implies that
  \[
    M^t_\infty (F_n \cap U_{n-1}) \geq M_\infty^t (U_{n-1}) > M_\infty^t (U) - \varepsilon.
  \]
  The set $(F_n \cap U_{n-1})_{(-\delta_n)}$ is open and increases to $F_n \cap U_{n-1}$
  as $\delta_n$ vanishes. Using Theorem~52 from \cite{rogers}, there is a $\delta_n$ such that
  \[
    M^t_\infty(U_n) = M_\infty^t ( (F_n \cap U_{n-1})_{(-\delta_n)} ) > M_\infty^t (U) - \varepsilon.
  \]

  We have that $\overline{U}_k \subseteq F_k$. Let $(C_l)$ be a cover by cylinders of
  $\cap \overline{U}_k$. Since $(\overline{U}_k)$ is a nested sequence of compact sets 
  and $C_l$ is open, there is an $m$ such that
  $\overline{U}_m \subseteq \cup C_l$. 

  Hence
  \[
    \sum_l d (C_l)^t \geq M_\infty^t (\overline{U}_m)
                        \geq M_\infty^t (U) - \varepsilon.
  \]
  As $\varepsilon$ is arbitrary we conclude that
  \[
    M_\infty^t \Bigl( \bigcap_i F_i \cap U \Bigr) \geq M_\infty^t (U).
  \]

  We now note that any finite intersection of open sets in $\mathcal G^t (\Lambda_A)$ is in $\mathcal G^t (\Lambda_A)$ and that any countable intersection of $G_\delta$ sets can be expressed as the intersection of a countable decreasing sequence of open sets.
  Together with the result proved above, this finishes the proof of the lemma.
\end{proof}

\section{Proof of Theorem~\ref{the:maintheorem}} \label{sec:proof}

Let $G_g (p)$ be the set of points in $\Sigma_A$ for which $p$ is an accumulation point of the ergodic averages of $g$, see (\ref{average}), and let $G_{g} (p, n, \varepsilon)$ be defined as
\[
  G_{g} (p, n, \varepsilon) = \Bigl\{\, \boldsymbol i \in \Sigma_A : \ p - \varepsilon < \frac{1}{n} \sum_{k = 0}^{n - 1} g(\sigma^k \boldsymbol i) < p + \varepsilon \, \Bigr\}.
\]
We also introduce the notation $\hat{G}_g (p)$ for the set of points for which the ergodic averages of $g$ converges to $p$. Hence $\hat{G}_g (p)$ is a subset of $G_g (p)$ and so $\dimH G_g (p) \geq \dimH \hat{G}_p (p)$.

The key step in proving Theorem \ref{the:maintheorem} and Corollary~\ref{cor:fulldim} is the following proposition.

\begin{proposition}\label{pro:inclass}
  For any real valued continuous function $g$ on $\Sigma_A$, it holds that $G_g (p) \in \mathcal G^t (\Sigma_A)$ for all $t < \dimH ( G_g(p) )$.
\end{proposition}

\begin{proof}[Proof of Theorem~\ref{the:maintheorem}]
  According to Proposition~\ref{pro:inclass}, the set $G_{g_i} (x_i) = \{ \, y \in \Lambda_A : x_i \in A_{g_i}(y) \, \}$ is in $\mathcal{G}^t (\Sigma_A)$ for $t < \dimH G_{g_i} (x_i)$. Now Theorem~\ref{the:maintheorem} follows from Theorem~\ref{the:classes}.
\end{proof}

\begin{proof}[Proof of Corollary~\ref{cor:fulldim}]
  Let $\varepsilon > 0$. 
  We will use Theorem~9 of \cite{barreirasaussol}. This theorem implies that there is a non-empty open interval $I_{g_i}$ such that $p \mapsto \dimH \hat{G}_{g_i} (p)$ is a real analytic function.

Since there are only finitely many maps $f_i$ and $f_{i,j}$ the conformal graph directed markov system is regular and there exists an ergodic measure of full dimension, see \cite{Mauldin-Urbanski2}.
By Birkhoffs ergodic theorem, it is therefore clear that $p \mapsto \dimH \hat{G}_{g_i} (p)$ attains the value $\dimH \Lambda_A$, somewhere in the interval $I_{g_i}$.
Hence, for each $g_i$ there are $x_{1,i} \ne x_{2,i}$ such that $\dimH \hat{G}_{g_i} (x_{1,i}) > \dimH \Lambda_A - \varepsilon$ and $\dimH \hat{G}_{g_i} (x_{2,i}) > \dimH \Lambda_A - \varepsilon$.
We then have $\dimH G_{g_i} (x_{1,i}) > \dimH \Lambda_A - \varepsilon$ and $\dimH G_{g_i} (x_{2,i}) > \dimH \Lambda_A - \varepsilon$.

Proposition~\ref{pro:inclass} implies that $G_{g_i} (x_{1,i})$ and $G_{g_i} (x_{2,i})$ are in $\mathcal{G}^{\dimH \Lambda_A - \varepsilon} (\Sigma_A)$.
Thereby $G_{g_i} (x_{1,i}) \cap G_{g_i} (x_{2,i})$ is in $\mathcal{G}^{\dimH \Lambda_A - \varepsilon} (\Sigma_A)$. This set consists of points for which both $x_{1,i}$ and $x_{2,i}$ are accumulation points of the ergodic averages of $g_i$.

Hence the set of points where the ergodic averages do not converge for any $g_i$ is in $\mathcal{G}^{\dimH \Lambda_A - \varepsilon} (\Sigma_A)$. Let $\varepsilon \to 0$.
\end{proof}

The proof of Proposition \ref{pro:inclass} is broken down to a series of lemmata.
The following lemma will be used to show that $G_g (p)$ is in the class $\mathcal G^t (\Sigma_A)$.

\begin{lemma} \label{lem:ings} 
  Let $(F_k)_{k=1}^\infty $ be a sequence of open sets in $\Sigma_A$. If there is a $c > 0$ such that
  \[
    \limsup_{k \to \infty} M^t_\infty (F_k \cap C) \geq c M_\infty^t (C)
  \]
  for all cylinders $C$ of generation $jm$, where $j \in \N$, then
  \[
    \bigcap_{n = 1}^\infty \bigcup_{k = n}^\infty F_k \in \mathcal{G}^{t'} (\Sigma_A), \  \forall t' \leq t.
  \]
\end{lemma}

\begin{proof}
  It follows from Lemma \ref{lem:removeconstant} that $\bigcup_{k = n}^\infty F_k \in \mathcal{G}^{t'} (\Sigma_A)$ for all $n$. The Lemma then follows from Proposition~\ref {pro:intersectionproperty}.
\end{proof}

We need to prove that the condition of Lemma~\ref{lem:ings} is satisfied for $G_g (p)$.
Of course, the sets $F_k$ in Lemma~\ref{lem:ings} will be the sets $G_g (p, k, \varepsilon)$. Hence, we will need that the sets $G_g (p, k, \varepsilon)$ are open. This is guarantied if $g$ is continuous. 
From now on we assume without mentioning that $g$ is a continuous function on $\Sigma_A$.

Before we continue we fix an $s < \dimH \Lambda_A$, and an $m = m(s)$ according to \eqref{eq:positive}. 
We define new outer measures $N_\infty^{m,t}$, $0 < t < s$, defined as $M_\infty^t$, but instead of considering all covers by cylinders, we only consider covers by cylinders of generation $km$, where $k\in \N$. We note that (\ref{eq:positive}) ensures that $N_\infty^{m,t} (C) = d (C)^t$ for all cylinders $C$ of generation $km$, where $k \in \N$. This makes it easier to work with $N_\infty^{m,t}$, than with $M_\infty^t$.

One observes that the measures $M_\infty^t$ and $N_\infty^{m,t}$ are equivalent. Hence there is some constant $c_1$, depending on $m$ and $t$, such that $M_\infty^t (U) \leq N_\infty^{m,t} (U) \leq c_1 M_\infty^t (U)$ holds for any open set $U$. Therefore, by Lemma~\ref{lem:removeconstant}, to prove that $G_g (p)$ is in the class $\mathcal G^t (\Sigma_A)$, it is sufficient to consider the measures $N_\infty^{m,t}$ instead of $M_\infty^t$.

\begin{lemma} \label{lem:Qa}
  Assume that there is a $c < 1$ and an integer $M$ such that
  \[
    N_\infty^{m,t} (C_{i_1 \dots i_m} \cap G_{g} (p, M, \varepsilon) ) < c N_\infty^{m,t} (C_{i_1 \dots i_m}),
  \]
  holds for the cylinder $C_{i_1 \dots i_m}$.
  Then for each cylinder $C = C_{c_1 \dots c_{jm}}$ such that $j\in \N$
  and $c_{(j-1)m+1}\dots c_{jm} ={i_1 \dots i_m}$, there is a constant $K_C$,
  depending only on $C$ and $g$, such that
  \[
    N_\infty^{m,t} (C \cap G_{g} (p, M, \varepsilon / 2)) < \kappa c N_\infty^{m,t} (C)
  \]
  holds if $M \geq K_C$.
\end{lemma}

\begin{proof}
  Let $C$ and $C_{i_1 \ldots i_m}$ be as in the statement of the lemma. Take $A$ such that $A > |g|$.

  By assumption, there is a cover of $C_{i_1 \dots i_m} \cap G_{g} (p, M, \varepsilon)$ by cylinders $(U_i)$ of generation $km$, with value less than $c N_\infty^{m,t} (C_{i_1 \dots i_m})$. For each $U_i$ there is a corresponding set $\tilde{U}_i$, created by appending $c_1 \ldots c_{(j-1)m}$ to the beginning of the coding of each element in $U_i$. 

  We note that
  $d (\tilde{U_i})^t \leq \kappa \frac{d (C)^t d (U_i)^t}{d (C_{i_1 \dots i_m})^t} = \kappa \frac{d (C)^t d (U_i)^t}{N_\infty^{m,t} (C_{i_1 \dots i_m})}$, by \eqref{eq:distortion2}. Moreover, $(\tilde{U}_i)_{i=1}^\infty$ is a cover of 
  \[
    C \cap G_{g} (p, M, \varepsilon - 2 A(j-1)m / M ).
  \]
  It now follows that
  \begin{multline*}
    N_\infty^{m,t} (C \cap G_{g} (p, M, \varepsilon - 2 A(j-1)m / M ) \leq \sum_i d ( \tilde{U}_i )^t \\
    = \frac{\kappa d ( C )^t }{N_\infty^{m,t} (C_{i_1 \dots i_m})}\sum_i d ( U_i )^t  < \kappa c d (C)^t = \kappa c N_\infty^{m,t} (C).
  \end{multline*}
  The lemma now follows if we choose $K_C$ so large that $2 A(j-1)m / K_C < \varepsilon / 2$.
\end{proof}

\begin{lemma} \label{lem:Qa2}
  Assume that there is a $c > 0$ and an integer $M$ such that
  \[
    N_\infty^{m,t} (C_{i_1 \dots i_m} \cap G_{g} (p, M, \varepsilon) ) > c N_\infty^{m,t} (C_{i_1 \dots i_m}),
  \]
  holds for the cylinder $C_{i_1 \dots i_m}$. Then for each cylinder $C=C_{c_1 \dots c_{jm}}$ such that $j\in \N$ and $c_{(j-1)m+1}\dots c_{jm} ={i_1 \dots i_m}$, there is a constant $K_C$, depending only on $C$ and $g$, such that
  \[
    N_\infty^{m,t} (C \cap G_{g} (p, M, 2\varepsilon )) > \kappa^{-1} c N_\infty^{m,t} (C),
  \]
  holds if $M \geq K_C$.
\end{lemma}

\begin{proof}
  Let $C$ and $C_{i_1 \ldots i_m}$ be as in the statement of the theorem.

  Let $A > |g|$ and consider the set $C_{i_1 \dots i_m} \cap G_{g} (p, M, \varepsilon)$. We can write this set as a union of cylinders. If we append $c_1 \ldots c_{(j-1)m}$ to the beginning of the coding of each of these cylinders, we get a collection of cylinders in 
  \[
    C \cap G_{g} (p, M, \varepsilon + 2 A (j-1)m / M).
  \]
  Let $(\tilde{U}_i)$ be a cover of $C \cap G_{g} (p, M, \varepsilon + 2 A (j-1)m / M)$. For each $\tilde{U}_i$ there is a corresponding set $U_i$, created by removing the first $(j-1)m$ symbols in the coding. 

  We note that $d (\tilde{U_i})t \geq \kappa^{-1} \frac{d (C)^t d (U_i)^t}{d (C_{i_1 \dots i_m})^t} = \kappa^{-1} \frac{d (C)^t d (U_i)^t}{N_\infty^{m,t} (C_{i_1 \dots i_m})}$ and that $(U_i)_{i=1}^\infty$ is a cover of $C_{i_1 \dots i_m} \cap G_{g} (p, M, \varepsilon)$.

  It now follows that
  \[
    \sum_i d ( \tilde{U}_i )^t = \frac{\kappa^{-1} d ( C )^t}{N_\infty^{m,t} (C_{i_1 \dots i_m})} \sum_i d ( U_i )^t
    > \kappa^{-1} c d (C)^t = \kappa^{-1} c N_\infty^{m,t} (C).
  \]
  Since the cover $(\tilde{U}_i)$ was arbitrary, we get
  \[
    N_\infty^{m,t} (C \cap G_{g} (p, M, \varepsilon + 2 A(j-1)m / M )) > \kappa^{-1} c N_\infty^{m,t} (C).
  \]
  The lemma now follows if we choose $K_C$ so large that $2 A(j-1)m / K_C < \varepsilon $.
\end{proof}

\begin{lemma}\label{lem:atob}
  There is a constant $L > 0$ such that if there is a $c > 0$ and an $M$ such that
  \[
    N_\infty^{m,t} (C_{i_1 \dots i_m} \cap G_{g} (p, M, \varepsilon)) > c  N_\infty^{m,t} (C_{i_1 \dots i_m}),
  \]
  for one cylinder $C_{i_1 \dots i_{m}}$ of generation $m$, then for each cylinder $C$, of generation $m$, there is a number $K_{C}$, depending only on $C$ and $g$, such that
  \[
    N_\infty^{m,t} (C \cap G_{g} (p, M, 2 \varepsilon)) > Lc N_\infty^{m,t} (C),
  \]
  if $M \geq K_C$.
\end{lemma}

\begin{proof}
  Since $\Sigma_A$ is transitive, there is a number $N \in \N$ such that any cylinder $C = C_{c_1 \dots c_m}$ contains a cylinder $C_{c_1 \dots c_{Nm}}$ of generation $Nm$ for which $c_{(N-1)m+1} \dots c_{Nm}=i_1 \dots i_m$. Since $\Sigma_A$ is of finite type and the maps $f = f_i$ or $f = f_{i,j}$ satisfy $\lambda_1 \leq |f'| \leq \lambda_2$, it is clear that there is a uniform constant $L' > 0$ such that
  \[
    \frac{d (C_{c_1 \dots c_{Nm}})^t}{d (C_{c_1 \dots c_m})^t} > L',
  \]
  regardless of which cylinders $C_{i_1 \dots i_{m}}$ and $C_{c_1 \dots c_{Nm}}$ we started with. (For instance, $L' = \lambda_1^{t (N - 1) m}$ will do.) Using Lemma \ref{lem:Qa2} we can find $K_C$ such that 
  \begin{align*}
    N_\infty^{m,t} (C_{c_1 \dots c_m} \cap G_{g} (p, M, 2 \varepsilon)) & \geq N_\infty^{m,t} (C_{c_1 \dots c_{Nm}} \cap G_{g} (p, M, 2 \varepsilon)) \\
    & \geq \kappa^{-1} c N_\infty^{m,t} (C_{c_1 \dots c_{Nm}}) > L' \kappa^{-1} c N_\infty^{m,t} (C_{c_1 \dots c_m})
  \end{align*}
  if $M \geq K_C$. Let $L = L' \kappa^{-1}$.
\end{proof}

\begin{lemma} \label{lem:stor}
  For each $t < \dimH (G_g (p))$ there is a constant $c > 0$ such that for any cylinder $C_{i_1 \dots i_m}$ of generation $m$ and any $\varepsilon > 0$, it holds that
  \[
    \liminf_{M \to \infty} N_\infty^{m,t} (C_{i_1 \dots i_m} \cap G_g (p, M, \varepsilon)) \geq c N_\infty^{m,t} (C_{i_1 \dots i_m}),
  \]
  for any $\varepsilon > 0$.
\end{lemma}

\begin{proof}
  Let $L$ be given by Lemma~\ref{lem:atob}.
  Assume on the contrary that there exist a cylinder $C_{i_1 \ldots i_m}$ and a strictly increasing sequence $(M_k)_{k=1}^\infty$ such that
  \begin{equation} \label{eq:contraryassumption}
    N_\infty^{m,t} (C_{i_1 \dots i_m} \cap G_{g} (p, M_k, 2 \varepsilon) )
    <  \kappa^{-1} \frac{L}{2} N_\infty^{m,t} (C_{i_1 \dots i_m}), \quad \forall k.
  \end{equation}
  By Lemma \ref{lem:atob} we get
  \[
    N_\infty^{m,t} ( C_{c_1 \dots c_m} \cap G_{g} (p, M_k, \varepsilon) )
    \leq \kappa^{-1} \frac{1}{2} N_\infty^{m,t} (C_{c_1 \dots c_m}), \quad \forall M_k \geq K_{C_{c_1 \dots c_m}}
  \]
  for each $C_{c_1 \dots c_m}$. Indeed, if would have
  \[
    N_\infty^{m,t} ( C_{c_1 \dots c_m} \cap G_{g} (p, M_k, \varepsilon) )
    > \kappa^{-1} \frac{1}{2} N_\infty^{m,t} (C_{c_1 \dots c_m}),
  \]
  for some $C_{c_1 \ldots c_m}$ and $M_k \geq K_{C_{c_1 \dots c_m}}$, then Lemma~\ref{lem:atob} gives us a contradiction to \eqref{eq:contraryassumption}.

  By Lemma \ref{lem:Qa} we get that
  \[
    N_\infty^{m,t} (C \cap G_p (M_{k}, \varepsilon/2 ) )
    \leq \frac{1}{2} N_\infty^{m,t} (C),\  \forall k \geq \max \{ K_C, K_{C_{c_1 \ldots c_m}} \}
  \]
  for each cylinder $C$ of any generation $jm$, $j \in \N$. Thus, there is a number $k_1$ and finite cover $(C_i)_i$ of $C_{i_1 \dots i_m} \cap G_g (p, M_{k_1}, \varepsilon / 2)$ such that
  \[
    \sum_i d (C_i)^t \leq \frac{2}{3} N_\infty^{m,t} (C_{i_1 \dots i_m}).
  \]

  Consider each $C_i$ and its intersection with $G_g (p, M_k, \varepsilon)$ for $k > k_1$.
  There is now a $k_2 > k_1$ and a finite cover $(C_{i,j})$ such that $(C_{i,j})_{j=1}^\infty$ covers $C_i \cap G_g (p, M_{k_2}, \varepsilon)$ and
  \[
    \sum_j d (C_{i,j})^t \leq \frac{2}{3} N_\infty^{m,t} (C_i) = \frac{2}{3} d (C_i)^t.
  \]
  We get $\sum_{i,j} d (C_{i,j})^t \leq \bigl( \frac{2}{3} \bigr)^2 N_\infty^{m,t} (C_{i_1 \dots i_m})$. 

  Continuing like this we get $\dimH (G_g(p)) < t$ which contradicts the assumptions.
\end{proof}

\begin{lemma}\label{lem:rightvalue}
For each $t < \dimH(G_g(p))$ there is a constant $c > 0$ such that for any cylinder $C_{i_1 \dots i_{jm}}$ of generation $jm$, where $j\in \N$, and any $\varepsilon > 0$, it holds that
\[
  \liminf_{M \to \infty} N_\infty^{m,t} (C_{i_1 \dots i_{jm}} \cap G_g (p, M, \varepsilon)) \geq c N_\infty^{m,t} (C_{i_1 \dots i_{jm}}).
\]
\end{lemma}

\begin{proof}
  This follows from Lemma \ref{lem:stor} and Lemma \ref{lem:Qa2}.
\end{proof}
  
We are now ready to prove Proposition \ref{pro:inclass}.

\begin{proof}[Proof of Proposition \ref{pro:inclass}]
  We first note that
  \[
    G_g (p) = \bigcap_{\varepsilon > 0} \bigcap_{n=1}^\infty \bigcup_{M=n}^\infty G_g (p, M, \varepsilon).
  \]
  By Theorem~\ref{the:classes}, Lemma \ref{lem:rightvalue} and Lemma \ref{lem:ings} it follows that $G_g (p)$ is in $\mathcal{G}^t (\Sigma_A)$ for any $t < \dimH (G_g (p))$.
\end{proof}

\section{Proof of Theorem~\ref{the:diophantine}} \label{sec:diophantineproof}

In this section we will prove Theorem~\ref{the:diophantine}, that $W (\alpha) \cap K$ is in $\mathcal{G}^t (K)$ for all $t < \frac{1}{\alpha} \frac{\log 2}{\log 3}$. We emphasise that $W (\alpha) \cap K$ is not in $\mathcal{G}^t (K)$ for any $t > \frac{1}{\alpha} \frac{\log 2}{\log 3}$, since $\dimH (W(\alpha) \cap K) \leq \frac{1}{\alpha} \frac{\log 2}{\log 3}$, as is easily shown by standard arguments.

Let
\[
  E_n = \{\, x \in [0,1) : |x - p/3^n| \leq 3^{-\alpha n}, \text{ for some } p \in \mathbb{Z} \,\},
\]
and let $F_n$ be the $n$-th step in the construction of the middle third Cantor set.
Then $E_n$ consists of $3^n + 1$ intervals of length $2 / 3^{\alpha n}$ (apart from two, which are of length $1 / 3^{\alpha n}$) and $F_n$ consists of $2^n$ intervals of length $1 / 3^n$. For any interval from $E_n$ that intersects $K$, the midpoints of the interval is an endpoint of an interval in $F_n$. Moreover, any endpoint of an interval in $F_n$ is also a midpoint of an interval in $E_n$. It follows that exactly $2^{n+1}$ intervals from $E_n$ intersect $K$ and these intersections are congruent.

Let $C$ be a cylinder in $K$ of generation $m$. Then $C$ is an interval of length $1 / 3^m$. We need to estimate
$\mathfrak{M}_\infty^t (C \cap K \cap W (\alpha))$. For this purpose we consider a cover of $C \cap K \cap W (\alpha)$ by cylinders $\{U_i\}$ in $K$. (This means that each $U_i$ is a projection by $\pi$ from a cylinder in the symbolic space $\Sigma_A$.) We may assume that this cover consists of sets with pairwise disjoint interior.

Consider one cylinder $U_i$ of generation $n_i$. We consider two cases. Either $n_i \leq n$ or $n_i > n$. In both cases we have $m$ smaller than $n$ and $n_i$.

If $n_i \leq n$ then $U_i$ intersects $2^{n - n_i + 1}$ intervals from $E_n$.

If $n_i > n$, then $U_i$ intersects at most one interval from $E_n$. We may then assume that $U_i$ intersect exactly one interval and that $U_i \leq 3 / 3^{\alpha n}$.

Let $\mu$ be a uniform mass distribution of mass 1 on $C \cap K \cap E_n$. In the first case when $n_i \leq n$ we have
\[
  \mu (U_i) = \frac{2^{n - n_i + 1}}{2^{n - m}} = 2 \frac{2^{-n_i}}{2^{-m}} \leq 2 \frac{|U_i|^t}{|C|^t},
\]
if $t \leq \frac{\log 2}{\log 3}$.

In case $n_i > n$ we first observe that if $I$ is the interval from $E_n$ that intersect $U_i$ then $\mu (U_i) / \mu (I) \leq (|U_i| / |I|)^{\frac{\log 2}{\log 3}}$. Hence
\begin{align*}
  \mu (U_i) & \leq \Bigl( \frac{|U_i|}{3^{-\alpha n}} \Bigr)^{\frac{\log 2}{\log 3}} \frac{1}{2^{n-m}}
  = \Bigl( \frac{|U_i|}{3^{-\alpha n}} \Bigr)^{\frac{\log 2}{\log 3}} \frac{1}{2^{n}} \frac{1}{|C|^\frac{\log 2}{\log 3}} \\
  &= \frac{|U_i|^t}{|C|^t} \Bigl( \frac{|U_i|}{3^{-\alpha n}} \Bigr)^{\frac{\log 2}{\log 3} - t} \frac{3^{\alpha n t}}{2^n} |C|^{t - \frac{\log 2}{\log 3}} \leq
 \frac{|U_i|^t}{|C|^t} \frac{3^{\alpha n t}}{2^n} |C|^{t - \frac{\log 2}{\log 3}}.
\end{align*}
Hence if $t < \frac{1}{\alpha} \frac{\log 2}{\log 3}$ and $n$ is sufficiently large, we have
\[
  \frac{3^{\alpha n t}}{2^n} |C|^{t - \frac{\log 2}{\log 3}} \leq 2,
\]
and so
\[
  \mu (U_i) \leq 2 \frac{|U_i|^t}{|C|^t}.
\]

Summing over $i$ we get
\[
  1 = \sum_i \mu (U_i) \leq \sum_i 2 \frac{|U_i|^t}{|C|^t},
\]
and $\sum_i |U_i|^t \geq \frac{1}{2} |C|^t$. Since $U_i$ is an arbitrary cover we get that $\mathfrak{M}_\infty^t (C \cap K \cap E_n) \geq \frac{1}{2} |C|^t$, provided that $t < \frac{1}{\alpha} \frac{\log 2}{\log 3}$ and $n$ is large. Now Lemma~\ref{lem:ings} implies that $\limsup K \cap E_n$ is in $\mathcal{G}^t (K)$ for $t < \frac{1}{\alpha} \frac{\log 2}{\log 3}$. This proves Theorem~\ref{the:diophantine}.

\end{document}